\documentclass[a4paper,12pt]{article}
  
\usepackage{amsfonts}
\usepackage{amssymb,amsmath} 
\numberwithin{equation}{section}
\newtheorem{theorem}{Theorem}[section]
\newtheorem{lemma}{Lemma}[section]

\newtheorem{definition}{Definition}[section]
\usepackage{cite}
\usepackage{hyperref}
\hypersetup{urlcolor=blue, colorlinks=true, citecolor=blue, linkcolor=blue}
\begin{document}
	\date{}
	\title{\textbf{On the Reversible Geodesics for a Finsler space with Randers change of Quartic metric}}
	\maketitle
	\begin{center}
		\author{\textbf{ Gauree Shanker, Ruchi Kaushik Sharma}}
	\end{center}
	
\textbf{Abstract.} In this paper, we consider a Finsler space with a Randers change of Quartic metric F = $\sqrt[4]{\alpha^4 + \beta^4} + \beta$. The conditions for this space to be with reversible geodesics are obtained. Further, we study some geometrical properties of F with reversible geodesics and prove that the Finsler metric F induces a generalized weighted quasi-distance $d_F$ on M.\\
\\
\textbf{M. S. C. 2010}: 58B20, 53B21, 53B40, 53C60.\\
{\bf Keywords -} Riemannian spaces; Reversible Geodesics; weighted quasi metric; Randers change; Quartic metric.
\section{Introduction}
An interesting topic in Finsler geometry is to study the reversible geodesics of a Finsler metric. Recall that, a Finsler space is said to have reversible geodesics if for every one of its oriented geodesic paths, the same path traversed in the opposite sense is also a geodesic. In the last decade many interesting and applicable  results have been obtained on the theory of Finsler spaces with reversible geodesics. In \cite{Cram05}, Crampin discussed Randers space with reversible geodesics. In (\cite{Shimada10}, \cite{Shimada13}), Masca, Sabau and Shimada have discussed reversible geodesics with ($\alpha$, $\beta$)-metric and two dimensional Finsler space with ($\alpha$, $\beta$)-metric to be with reversible geodesic, respectively. In \cite{Shimada12}, Sabau and Shimada have given some important results on reversible geodesics. In \cite{Shnaker16}, Shanker and Baby have discussed reversible geodesics for generalized ($\alpha, \beta$)-metric. Recently, Shanker and Rani \cite{Shanker18} have studied weighted quasi metric associated with Finsler spaces with reversible geodesics.\\
 In this paper, we find conditions for a Finsler space (M, F) with Randers change of Quartic metric F = $\sqrt[4]{\alpha^4 + \beta^4} + \beta$ to be with reversible geodesics. The main results of this paper lies in theorem (3.1), (4.1), (4.2), (5.1) and (5.2). \\
\section{Preliminaries}
Let $F^n$ = ($M, F$) be a connected \textit{n}-dimensional Finsler manifold and let TM = $\bigcup\limits_{x \in M} T_xM$ denotes the tangent bundle of M with local coordinates u = (x, y) = $(x^i, y^i) \in TM$, where \textit{i} = 1, ..., \textit{n}, \textit{y} = $y^i\dfrac{\partial}{\partial x^i}$.
\par 
If $\gamma$ : [0, 1] $\longrightarrow$ M is a piecewise $C^\infty$ curve on M, then its Finslerian length is defined as  
\begin{align}
L_F(\gamma) = \int_{0}^{1} F(\gamma(t), \dot{\gamma}(t))dt,
\end{align}
and the Finslerian distance function $d_F$ : M $\times$ M $\longrightarrow$ [0, $\infty$) is defined by $d_F(p,q) = inf_\gamma L$, where infimum is taken over all piecewise $C^\infty$ curves $\gamma$ on M joining the points p, q $\in$ M.
In general, this is not symmetric.
\par 
A curve $\gamma$ : [0, 1] $\longrightarrow$ M is called a geodesic of (M, F) if it minimizes the Finslerian length for all piecewise $C^\infty$ curves that keep their endpoints fixed. We denote the reverse Finsler metric of F as $\tilde{F}$ : TM $\longrightarrow$ (0, $\infty$), given by $\tilde{F}(x, y)$ = $F(x, -y)$. One can easily see that $\tilde{F}$ is also a Finsler metric.
\\
\begin{lemma}
A Finsler metric is with a reversible geodesic if and only if for any geodesic $\gamma$(t) of F, the reverse curve $\tilde{\gamma}(t)$ = $\gamma$(1-\textit{t}) is also a geodesic of F.	
\end{lemma}
\begin{lemma}
	Let $(M, F)$ be a connected, complete Finsler manifold with associated distance function $d_F$ : \textit{M} $\times$ \textit{M} $\longrightarrow$ [0, $\infty$). Then, $d_F$ is a symmetric distance function on \textit{M} $\times$ \textit{M} if and only if \textit{F} is a reversible Finsler metric, i.e., $F(x, y) = F(x, -y)$.
\end{lemma}
\begin{lemma}
	A smooth curve $\gamma$ : [0, 1] $\longrightarrow$ \textit{M} is a constant Finslerian speed geodesic of $(M, F)$ if and only if it satisfies $\ddot{\gamma}$ + 2 $G^i$($\gamma(t)$, $\dot{\gamma}(t)$) = 0, \textit{i} = 1, ..., \textit{n}, where the functions $G^i$ : \textit{TM} $\longrightarrow \mathbb{R}$ are given by 
	\begin{equation}
	G^i(x, y) = \Gamma_{jk}^{i}(x, y)y^iy^j,
	\end{equation}
	with $\Gamma_{jk}^{i}(x, y) = \dfrac{g^{is}}{2}\Bigl(\dfrac{\partial g_{sj}}{\partial x^k} + \dfrac{\partial g_{sk}}{\partial x^j} - \dfrac{\partial g_{jk}}{\partial x^s}\Bigr)$.
\end{lemma}
\textbf{Remark 1.} It is well known \cite{Shen01} that the vector field $\Gamma$ = $y^i\dfrac{\partial}{\partial x^i} - 2G^i \dfrac{\partial}{\partial y^i}$, is a vector field on TM, whose integral lines are the canonical lifts $\tilde{\gamma}(t)$ = ($\gamma$(t), $\dot{\gamma}(t)$) of the geodesics of $\gamma$. This vector field $\Gamma$ is called the canonical geodesics spray of the Finsler space $(M, F)$ and $G^i$ are called the coefficients of the geodesics spray $\Gamma$.
\begin{definition}
	If F and $\tilde{F}$ are two different fundamental Finsler functions on the same manifold M, then they are said to be projectively equivalent if their geodesics coincide as set points.
\end{definition}
\begin{lemma}
	A Finsler structure ($M, F$) is with a reversible geodesic if and only if \textit{F} and its reverse function $\tilde{F}$ are projectively equivalent. 
\end{lemma}
\section{Reversible Geodesics for a Finsler space with Randers change of Quartic metric.}
Consider a Finsler space ($M, F$) with a special ($\alpha, \beta$)-metric \textit{F} = $\sqrt[4]{\alpha^4 + \beta^4} + \beta$. Here, \textit{F} can be treated as the Randers change of Quartic-metric $\tilde{F} = \sqrt[4]{\alpha^4 + \beta^4}$. One can easily see that $\tilde{F}(x, -y) = \tilde{F}(x, y)$.\\
As we know that \cite{Shimada12} if $(M, F)$ is a non-Riemannian \textit{n}(\textit{n} $\geq$ 2)-dimensional Finsler space with ($\alpha, \beta$)-metric, which is not absolute homogeneous, then \textit{F} is with reversible geodesics if and only if $F(\alpha, \beta) = F_0(\alpha, \beta) + \epsilon \beta$, where $F_0$ is absolute homogeneous ($\alpha, \beta$)-metric, $\epsilon$ is a non-zero constant and $\beta$ is a closed 1-form on the Manifold \textit{M}.
\par 
In our case, $F_0$ = $\tilde{F}$, which is absolute homogeneous. If $\beta$ is a closed 1-form, then \textit{F} is with reversible geodesics. Further, a necessary and sufficient condition for \textit{F} to have reversible geodesics is that \cite{Shimada12}
\begin{align}
\tilde{\Gamma}\Bigl(\dfrac{\partial F}{\partial y^i}\Bigr) - \dfrac{\partial F}{\partial x^i} = 0,
\end{align}
where $\tilde{\Gamma}$ is the reverse of $\Gamma$, the geodesic spray of \textit{F}; moreover $\tilde{\Gamma}$ is geodesic spray of $\tilde{F}$. We have F = $\tilde{F}$ + $\beta$, where,
 $\tilde{F} = \sqrt[4]{\alpha^4 + \beta^4}$.

Therefore,  
\begin{align}
&\tilde{\Gamma}\Bigl(\dfrac{\partial F}{\partial y^i}\Bigr) - \dfrac{\partial F}{\partial x^i} \nonumber \\
 &= \tilde{\Gamma}\Bigl(F_{\alpha} \dfrac{\partial \alpha}{\partial y^i} + F_{\beta} \dfrac{\partial \beta}{\partial y^i}\Bigr) - F_{\alpha} \dfrac{\partial \alpha}{\partial x^i} - F_{\beta}\dfrac{\partial \beta}{\partial x^i} \nonumber \\
&
= \tilde{\Gamma}(F_{\alpha})\dfrac{\partial \alpha}{\partial y^i} + F_{\alpha} \tilde{\Gamma}\Bigl(\dfrac{\partial \alpha}{\partial y^i}\Bigr) + \tilde{\Gamma}(F_{\beta})\Bigl(\dfrac{\partial \beta}{\partial y^i}\Bigr) + F_{\beta}\tilde{\Gamma}\Bigl(\dfrac{\partial \beta}{\partial y^i}\Bigr) - F_{\alpha} \dfrac{\partial \alpha}{\partial x^i} - F_{\beta} \dfrac{\partial \beta}{\partial x^i} \nonumber \\
&
=\tilde{\Gamma}(F_{\alpha}) \dfrac{\partial \alpha}{\partial y^i} + F_{\alpha} \Bigl[\tilde{\Gamma}\Bigl(\dfrac{\partial \alpha}{\partial y^i}\Bigr) - \dfrac{\partial \alpha}{\partial x^i}\Bigr] + \tilde{\Gamma}(F_{\beta}) \dfrac{\partial \beta}{\partial y^i} + F_{\beta}\Bigl[\tilde{\Gamma} \Bigl(\dfrac{\partial \beta}{\partial y^i}\Bigr) - \dfrac{\partial \beta}{\partial x^i}\Bigr].
\end{align}
For the Riemannian metric $\alpha$, the Euler-Lagrange equation gives $\tilde{\Gamma}\Bigl(\dfrac{\partial \alpha}{\partial y^i}\Bigr) - \dfrac{\partial \alpha}{\partial x^i}$ = 0. Also, one knows \cite{Shimada12} that if (\textit{M}, \textit{F}($\alpha, \beta$)) is a Finsler space with ($\alpha$, $\beta$)-metric, then $f(x, y)\dfrac{\partial \alpha}{\partial y^i} + g(x, y)b_i$ = 0, $\forall$ \textit{i} = 1, 2, ..., \textit{n}, implies that \textit{f} = \textit{g} = 0, for any smooth functions \textit{f} and \textit{g} on \textit{TM}. It is known that, if $\beta$ is closed and \textit{F} is projectively equivalent to the Riemannian metric $\alpha$, then  \\ $\tilde{\Gamma}(F_{\alpha})\dfrac{\partial \alpha }{\partial y^i} + \tilde{\Gamma}(F_{\beta})b_i$ = 0 
\\ and hence by using lemma 2.7 of \cite{Shimada12}, we find that $\tilde{\Gamma}(F_{\alpha}) = 0, \tilde{\Gamma}(F_{\beta}) = 0.$\\
Again, since \textit{F} = $\sqrt[4]{\alpha^4 + \beta^4} + \beta$, therefore $F_{\beta} = 1 + \dfrac{\beta^3}{(\alpha^4 + \beta^4)^{\frac{3}{4}}}$.\\
Now, using the above results, the equation (3.2) reduces to the form
\begin{eqnarray}
\tilde{\Gamma}\Bigl(\dfrac{\partial F}{\partial y^i}\Bigr) - \dfrac{\partial F}{\partial x^i} &=&  F_{\beta}\Bigl[\tilde{\Gamma} \dfrac{\partial \beta}{\partial y^i} - \dfrac{\partial \beta}{\partial x^i}\Bigr] \nonumber \\
&&
= \Bigl(1 + \dfrac{\beta^3}{(\alpha^4 + \beta^4)^{\frac{3}{4}}}\Bigr)\Bigl[\bar{\Gamma}(b_i) - \dfrac{\partial b_j}{\partial x^i}y^j\Bigr] \nonumber \\
&&
= \Bigl(1 + \dfrac{\beta^3}{(\alpha^4 + \beta^4)^{\frac{3}{4}}}\Bigr)\Bigl[\dfrac{\partial b_i}{\partial y^j}y^j - \dfrac{\partial b_j}{\partial x^i} y^j\Bigr] \nonumber \\
&&
= \Bigl(1 + \dfrac{\beta^3}{(\alpha^4 + \beta^4)^{\frac{3}{4}}}\Bigr)\Bigl[\dfrac{\partial b_i}{\partial x^j} - \dfrac{\partial b_j}{\partial x^i}\Bigr]y^j. 
\end{eqnarray}
Now, $1 + \dfrac{\beta^3}{(\alpha^4 + \beta^4)^{\frac{3}{4}}}$ can not be zero. Therefore, from equation (3.1) and (3.3) we conclude that
 F is with reversible geodesics 
 if and only if
$\Bigl(\dfrac{\partial b_i}{\partial x^j} - \dfrac{\partial b_j}{\partial x^i}\Bigr)y^j$ = 0.\\
i.e., $\tilde{F}$ is with reversible geodesic if and only if $\beta$ is closed 1-form. Hence, we have the following theorem:
\begin{theorem}
	A Finsler space $(M, F)$ with Randers change of Quartic metric F = $\tilde{F} + \beta, where, \tilde{F} = \sqrt[4]{(\alpha^4 + \beta^4)}$, is with reversible geodesics if and only if the differential 1-form $\beta$ is closed on \textit{M}.
\end{theorem}

\section{Projective Flatness of Randers change of Quartic metric}
A Finsler space $(M, F)$ is called (locally) projectively flat if all its geodesics are straight lines \cite{Shen09}. An equivalent condition is that the spray coefficients $G^i$ of \textit{F} can be expressed as $G^i = P(x, y)y^i$, where $P(x, y)$ = $\dfrac{1}{2F}\dfrac{\partial F}{\partial x^k}y^k$.
\par 
An equivalent characterization of projective flatness is the Hamel's relation \cite{Hamel03}
\begin{align*}
\dfrac{\partial^2 F}{\partial x^m \partial y^k}y^m - \dfrac{\partial F}{\partial x^k} = 0.
\end{align*}
Recall that (\cite{Shnaker16}, \cite{Shanker18}) if F = $F_0$ + $\epsilon$ $\beta$ is a Finsler metric, where $F_0$ is an absolute homogeneous ($\alpha$, $\beta$)-metric, then any two of the following properties imply the third one:\\
\\
(1) \textit{F} is projectively flat;\\
\\
(2) $F_0$ is projectively flat;\\
\\
(3) $\beta$ is closed.\\
In our case, F = $\tilde{F}$ + $\beta$, where $\tilde{F} = \sqrt[4]{(\alpha^4 + \beta^4)}$, which is absolute homogeneous. Hence, we have the following:
\begin{theorem}
	Let $(M, F)$ be a Finsler space with Randers change of Quartic metric \textit{F} = $\sqrt[4]{(\alpha^4 + \beta^4)} + \beta$. Then, \textit{F} is projectively flat if and only if $\tilde{F}$ is projectively flat. 
	\end{theorem}
\textit{Proof:} Let $(M, F)$ be projectively flat, then by Hamel's relation for projective flatness, we have 
\begin{align*}
\dfrac{\partial^2 F}{\partial x^m \partial y^k}y^m - \dfrac{\partial F}{\partial x^k} = 0.
\end{align*}
The proof of the theorem directly follows from it.
\begin{theorem}
Let $(M, F)$ be a Finsler space with Randers change of Quartic metric. If \textit{F} is projectively flat, then it is with reversible geodesics.
\end{theorem}
\textit{Proof.} Applying Hammel's equation, one can easily see that \textit{F} is projectively flat if and only if $\tilde{F}$ is projectively flat, which implies that both F and $\tilde{F}$ are projectively equivalent to the standard Euclidean metric and therefore F must be projective to $\tilde{F}$. Thus, \textit{F} must be with a reversible geodesic.
\section{Weighted quasi metric associated with Randers change of Quartic metric}
It is well known that the Riemannian spaces can be represented as metric spaces. Indeed, for a Riemannian space (\textit{M}, $\alpha$), one can define the induced metric space (\textbf{M}, $d_{\alpha}$) with the metric\\
\begin{eqnarray}
d_{\alpha} : M \times M \longrightarrow [0, \infty), d_{\alpha}(x, y) = \inf\limits _{\gamma \in \Gamma_{xy}} \int_{a}^{b} \alpha(\gamma(t), \dot{\gamma}(t))dt,
\end{eqnarray}
where $\Gamma_{xy}$ = \{$\gamma$ : [\textit{a}, \textit{b}] $\longrightarrow$ \textit{M} $\mid$ $\gamma$ is piecewise, $\gamma$(a) = \textit{x}, $\gamma$(b) = \textit{y}\} is the set of curves joining x and y, $\dot{\gamma}(t)$ is the tangent vector to $\gamma$ at $\gamma$(t). Then $d_{\alpha}$ is a metric on \textit{M} satisfying the following conditions:\\
\\
1. Positiveness : $d_{\alpha}(x, y)$ $>$ 0, if x $\neq$ y, $d_{\alpha}$(\textit{x}, \textit{x}) = 0, \textit{x}, \textit{y} $\in$ X. \\
\\
2. Symmetry : $d_{\alpha}$(\textit{x}, \textit{y}) = $d_{\alpha}$(\textit{y}, \textit{x}), $\forall$ \textit{x}, \textit{y} $\in$ \textit{M}.\\
\\
3. Triangle inequality: $d_{\alpha}$(\textit{x}, \textit{y}) $\leq$ $d_{\alpha}$(x, z) + $d_{\alpha}$(\textit{z}, \textit{y}), $\forall$ \textit{x}, \textit{y}, \textit{z} $\in$ \textit{M}. \\
\\
Similar to the Riemannian space, one can induce the metric $d_F$ to a Finsler space (\textit{M}, \textit{F}), given by
\begin{eqnarray}
d_F : M \times M \longrightarrow [0, \infty), d_F(x, y) = \inf\limits _{\gamma \in \Gamma_{xy}} \int_{a}^{b} F(\gamma(t), \dot{\gamma}(t))dt.
\end{eqnarray}
But unlike the Riemannian case, here $d_F$ lacks the symmetric condition. In fact, $d_F$ is a special case of quasi metric defined below:
\begin{definition}
	A quasi metric \textit{d} on a set \textit{X} is a function \textit{d} : \textit{X} $\times$ \textit{X} $\longrightarrow$ [0, $\infty$) that satisfies the following axioms:\\
	\\
	1. Positiveness : \textit{d}(\textit{x}, \textit{y}) $>$ 0, if \textit{x} $\neq$ \textit{y}, \textit{d}(\textit{\textit{x}}, \textit{x}) = 0, \textit{x}, \textit{y} $\in$ \textit{X}. \\
	\\
	2. Triangle inequality : \textit{d}(\textit{x}, \textit{y}) $\leq$ \textit{d}(\textit{x}, \textit{z}) + \textit{d}(\textit{z}, \textit{y}), $\forall$ \textit{x}, \textit{y}, \textit{z} $\in$ \textit{X}. \\
	\\
	3. Separation axiom : \textit{d}(\textit{x}, \textit{y}) = \textit{d}(\textit{y}, \textit{x}) = 0 $\Rightarrow$ \textit{x} = \textit{y}, $\forall$ \textit{x}, \textit{y} $\in$ \textit{X}.\\
	\\
\end{definition}
One special class of quasi metric spaces are the so called weighted quasi metric spaces (\textit{M, \textit{d}, \textit{w}}), where \textit{d} is a quasi-metric on \textit{M} and for each \textit{d}, there exist a function \textit{w} : \textit{M} $\longrightarrow$ [0, $\infty$), called the \textit{weight} of \textit{d} that satisfies
\\
4. \textit{Weightability : d(x, y) + w(x) = d(y, x) + w(y), $\forall$  x, y $\in$ M}. \par 
In this case, the weight function w is $\mathbb{R}$-valued, and is called generalized weight.
\begin{theorem}
	Let (\textit{M}, \textit{F}) be an \textit{n}-dimensional simply connected smooth Finsler manifold with \textit{F} as Randers change of Quartic metric. Then, \textit{F} induces generalized weighted quasi metric $d_F$ on \textit{M}.
\end{theorem}
\textit{Proof.} We consider that (\textit{M}, \textit{F}) is a Finsler space with F = $\beta + \sqrt[4]{(\alpha^4 + \beta^4)}$, which can be written as \textit{F} = $\tilde{F}$ + $\beta$, where $\tilde{F}$ = $\sqrt[4]{(\alpha^4 + \beta^4)}$ is an 
absolute homogeneous Finsler metric on M and $\beta$ an exact 1-form.
\par 
Let $\gamma_{xy}$ $\in$ $\Gamma_{xy}$ be an \textit{F}-geodesic, which is at the same time $\tilde{F}$-geodesic, then from equation (5.2), we get
\begin{eqnarray}
d_F(x, y) &=& \int_{a}^{b} F(\gamma(t), \dot{\gamma}(t))dt \nonumber \\
&&
= \int_{a}^{b} \Bigl( \beta + \sqrt[4]{(\alpha^4 + \beta^4)}\Bigr)dt \nonumber \\
&&
= \int_{a}^{b}\Bigl(\sqrt[4]{(\alpha^4 + \beta^4)}\Bigr)dt + \int_{a}^{b}\beta dt \nonumber \\
&&
= \int_{\gamma_{xy}}\Bigl(\sqrt[4]{(\alpha^4 + \beta^4)}\Bigr) + \int_{\gamma_{xy}} \beta. 
\end{eqnarray}
Consider a fixed point $a \in M$ and define the function $w_{a} : M \longrightarrow \mathbb{R}$ by $w_{a}(x) := d_{F}(a, x) - d_{F}(x, a)$.\\
From the equation (5.3) it follows that 
\begin{align}
w_{a}(x) = \int_{\gamma_{ax}} \beta - \int_{\gamma_{xa}} \beta = -2\int_{\gamma_{xa}} \beta,
\end{align}
where we have used the Stokes theorem for the 1-form $\beta$ on the closed domain D with boundary $\partial$D := $\gamma_{ax} \bigcup \gamma_{xa}$. \\
It can be easily seen that $w_{a}$ is an anti-derivative of $\beta$. This is well defined if and only if the integral in the R.H.S. of equation(5.4) is path independent, i.e., $\beta$  must be exact.
\par 
Then $d_F$ is a weighted quasi-metric with generalized weight $w_{a}$. Next we have
\begin{eqnarray}
d_F(x, y) + w_{a}(x) = \int_{\gamma_{xy}} \Bigl(\sqrt[4]{(\alpha^4 + \beta^4)}\Bigr) - \int_{\gamma_{xa}} \beta - \int_{\gamma_{ya}}\beta,
\end{eqnarray} 
where we have again used the Stokes theorem for the one form $\beta$ on the closed domain with boundary $\gamma_{ax} \bigcup \gamma_{xy} \bigcup \gamma_{ya}$.
\par 
Similarly,
\begin{eqnarray}
d_F(y, x) + w_{a}(y) = \int_{\gamma_{yx}} \Bigl(\sqrt[4]{(\alpha^4 + \beta^4)}\Bigr) - \int_{\gamma_{ya}} \beta - \int_{\gamma_{xa}}\beta.
\end{eqnarray}
 From the equations (5.5) and (5.6) we conclude that $d_F$ is weighted quasi-metric with generalized weight $w_a$. \\
 This completes the proof.
\\
Next, recall the following:
\begin{lemma} (\cite{Kunzi94}, \cite{Vitolo95}) 
	Let (\textit{M}, \textit{d}) be any quasi-metric space. Then \textit{d} is weightable if and only if there exists \textit{w} : \textit{M} $\longrightarrow$ [0, $\infty$) such that
	\begin{eqnarray}
	d(x, y) = \rho(x, y) + \dfrac{1}{2}[w(x) - w(y)], \forall x, y \in M,
	\end{eqnarray} 
	where $\rho$ is the symmetrized distance function of \textit{d}. Moreover, we have
	\begin{eqnarray}
	\dfrac{1}{2}[w(x) - w(y)] \leq \rho(x, y), \forall x, y \in M.
	\end{eqnarray}
\end{lemma}
The proof is trivial from the definition of weighted quasi-metric.\\
\textbf{Remark 2.} If (\textit{M}, \textit{F}) is a Finsler space with a special ($\alpha$, $\beta$)-metric \textit{F} = $\sqrt[4]{(\alpha^4 + \beta^4)} + \beta$, then the induced quasi-metric $d_F$ and the symmetrized metric $\rho$ induce the same topology on \textit{M}. This follows immediately from (\cite{Shen01}, \cite{Shen09}).\\
\textbf{Remark 3.} From Lemma 5.1, It can be seen that the assumption of \textit{w} to be smooth is not essential.\\
Next, we discuss an interesting geometric property concerning the geodesic triangles.
\begin{theorem}
	Let (\textit{M}, \textit{F}) be a Finsler space with the Randers change of Quartic-metric \textit{F} = $\sqrt[4]{(\alpha^4 + \beta^4)} + \beta$. Then the parameteric length of any geodesic triangle on \textit{M} does not depend on the orientation, that is,
	\begin{align*}
	d_F(x, y) + d_F(y, z) + d_F(z, x) = d_F(x, z) + d_F(z, y) + d_F(y, x), \forall x, y, z \in M. 
	\end{align*}
\end{theorem}
\textit{Proof.} Since the Randers change of Quartic metric F = $\sqrt[4]{(\alpha^4 + \beta^4)} + \beta$ can be treated as the Randers change of absolute homogeneous Finsler metric $\tilde{F} = \sqrt[4]{(\alpha^4 + \beta^4)}$, i.e., F = $\tilde{F}$ + $\beta$ with d$\beta$ = 0, from theorem 5.1, it follows that the quasi-metric is weightable and therefore equation (5.7) holds good. By using the formula (5.7), a simple calculation gives the required result.

\newpage
Authors' addresses:\\
\\
Gauree Shanker\\
Department of Mathematics and Statistics,\\
School of Basic and Applied Sciences,\\
Central University of Punjab, Bathinda-151 001, \\
Punjab, India.\\
E-mail:  grshnkr2007@gmail.com\\ 
\\
\\	
R. K. sharma, \\
Department of Mathematics and Statistics \\
Banasthali University, Banasthali-304022\\ 
Rajasthan, India.\\
E-mails: ruchikaushik07@gmail.com

\begin{thebibliography}{99}
	\bibitem{Cram05}
	M. Crampin, \emph{Randers spaces with reversible geodesics}, Publ. Math. Debrecen, \textbf{67/3}(2005), 401-409.
	\bibitem{Hamel03}
	G. Hamel, \emph{$\ddot{U}ber$ die Geomtrieen in denen die Geraden die K$\ddot{u}rzeinten $sind}, 
	Math. Ann., \textbf{57}(1903) .
	\bibitem{Kunzi94}
	H. P. A. Kunzi, V. Vajner, \emph{Weighted quasi-metrics},
	Annals New York Acad. Sci., \textbf{728}(1994), 64-77.
	\bibitem{Shimada10}
	I. Masca, S. V. Sabau, H. Shimada, \emph{Reversible geodesics for ($\alpha$, $\beta$)-metrics}, Intl. J. Math., \textbf{21}(2010), 1071-1094.
	\bibitem{Shimada13}
	I. Masca, S. V. Sabau, H. Shimada, \emph{Two dimensional ($\alpha$, $\beta$)-metrics with reversible geodesics}, Publ. Math. Debrecen, \textbf{82}(2013), 485-501.
	\bibitem{Shimada12} 
	S. V. Sabau, H. Shimada \emph{Finsler manifolds with reversible geodesics}, Rev. Roumaine Math. Pures Appl, \textbf{57}(2012), 91-103.
	\bibitem{Shen01} 
	Z. Shen, \emph{Differential geometry of sprays and Finsler spaces}, Kluwer Academic Publishers, Dordrecht, (2001). 
	\bibitem{Shen09}
	Z. Shen, \emph{On projectively flat ($\alpha$, $\beta$)-metrics}, Canad. Math. Bull.,
	\textbf{52}(2009), 132-144.
	\bibitem{Vitolo95} P. Vitolo, \emph{A representation theorem for quasi-metric space}, Topology Appl., \textbf{65}(1995), 101-104.
	\bibitem{Shnaker16} 
	G. Shanker, S. A. Baby, \emph{Reversible geodesics of Finsler spaces with a special ($\alpha$, $\beta$)-metric}, Bull. Cal. Math. Soc., \textbf{109}(2017), 183-188.
	\bibitem{Shanker18}
	G. Shanker, S. Rani, Weighted quasi-metrics associated with Finsler metrics, arXiv:1801.05636v2[math.DG]. 
\end{thebibliography}
\end{document}